\documentclass[a4paper]{scrartcl}\usepackage[]{graphicx}\usepackage[]{color}
\makeatletter
\def\maxwidth{ %
  \ifdim\Gin@nat@width>\linewidth
    \linewidth
  \else
    \Gin@nat@width
  \fi
}
\makeatother

\definecolor{fgcolor}{rgb}{0.345, 0.345, 0.345}

\usepackage{framed}
\makeatletter
 {\par\unskip\endMakeFramed%
 \at@end@of@kframe}
\makeatother

\definecolor{shadecolor}{rgb}{.97, .97, .97}
\definecolor{messagecolor}{rgb}{0, 0, 0}
\definecolor{warningcolor}{rgb}{1, 0, 1}
\definecolor{errorcolor}{rgb}{1, 0, 0}

\usepackage{alltt}
       
\usepackage[utf8]{inputenc}
\usepackage[T1]{fontenc}
\usepackage[english]{babel} 
\usepackage{lmodern}
\usepackage[]{algorithm2e}
\usepackage[onehalfspacing]{setspace}
\usepackage{amsmath,amssymb,amsthm,amsfonts,amsbsy,latexsym}
\usepackage{graphicx}
\usepackage{geometry}
\usepackage{natbib}
\usepackage{multirow}
\usepackage{float}
\usepackage{fancyhdr}
\usepackage{makecell}
\usepackage{pgf,tikz}
\usepackage{makeidx}
\usetikzlibrary{shapes,snakes}
\usetikzlibrary{fadings}
\usepackage[pdfpagelabels=true]{hyperref}

\newcommand{\E}{\mbox{I\negthinspace E}}
\newcommand{\B}{\mbox{I\negthinspace B}}

\newcommand{\R}{\mathbb{R}}

\newcommand{\N}{\mathbb{N}}

\newcommand{\D}{\mathcal{D}}

\DeclareMathOperator*{\argmin}{argmin}
\DeclareMathOperator*{\argmax}{argmax}

\DeclareMathOperator*{\IC}{IC}

\DeclareMathOperator*{\runs}{runs}
\DeclareMathOperator*{\mean}{mean}

\DeclareMathOperator*{\tr}{tr}

\providecommand{\keywords}[1]
{
  \small	
  \textbf{\textit{Keywords---}} #1
}

\numberwithin{equation}{section}
\IfFileExists{upquote.sty}{\usepackage{upquote}}{}
\begin{document}

\newtheorem{thm}{Theorem}[section]
\newtheorem{Def}[thm]{Definition}
\newtheorem{lem}[thm]{Lemma}
\newtheorem{rem}[thm]{Remark}
\newtheorem{cor}[thm]{Corollary}
\newtheorem{ex}[thm]{Example}
\newtheorem{ass}[thm]{Assumption}
\newtheorem*{bew}{Proof}
\newtheorem{prop}[thm]{Proposition}
\title{The column measure and Gradient-Free Gradient Boosting}
\author{Tino Werner\footnote{Carl von Ossietzky University Oldenburg}; Peter Ruckdeschel \footnote{Carl von Ossietzky University Oldenburg}}
\maketitle

\begin{footnotesize} 

\begin{abstract}

Sparse model selection by structural risk minimization leads to a set of a few predictors, ideally a subset of the true predictors. This selection clearly depends on the underlying loss function $\tilde L$. For linear regression with square loss, the particular (functional) Gradient Boosting variant $L_2-$Boosting excels for its computational efficiency even for very large predictor sets, while still providing suitable estimation consistency. For more general loss functions, functional gradients are not always easily accessible or, like in the case of continuous ranking, need not even exist. To close this gap, starting from column selection frequencies obtained from $L_2-$Boosting, we introduce a loss-dependent ''column measure'' $\nu^{(\tilde L)}$ which mathematically describes variable selection. The fact that certain variables relevant for a particular loss $\tilde L$ never get selected by $L_2-$Boosting is reflected by a respective singular part of $\nu^{(\tilde L)}$ w.r.t. $\nu^{(L_2)}$. With this concept at hand, it amounts to a suitable change of measure (accounting for singular parts) to make $L_2-$Boosting select variables according to a different loss $\tilde L$. As a consequence, this opens the bridge to applications of simulational techniques such as various resampling techniques, or rejection sampling, to achieve this change of measure in an algorithmic way. 

\end{abstract}

\end{footnotesize} \


\begin{small} \ \\ \

\section{Introduction} \

This work is motivated by an application in the context of tax auditing, where due to resource restrictions in terms of auditors, in a continuous ranking problem, one seeks to find the tax declaration prone to the highest level of tax evasion, giving rise to a highly non-differentiable loss function. Model selection from a set of candidate predictive criteria (usually categorical variables) is to respect this particular non-smooth loss function. Similar use cases arise in general fraud detection and resource allocation. For reference, see \cite{pickett} for a broad overview, \cite{moraru} for a short survey of different risks in auditing and \cite{khanna} and \cite{bowlin} for a study on bank-internal risk-based auditing resp. for a study on risk-based auditing for resource planning. \ \\

It has turned out in several works (\cite{alm}, \cite{gupta}, \cite{hsu15}) that risk-based auditing based on machine learning algorithms is far more sophisticated than simply selecting instances randomly. Since one is interested in an ordering of the instances, this setting corresponds to ranking problems.  \ \\

More specifically, in our tax auditing example we head for the hard continuous ranking problem, i.e., we want to order all observations according to a continuous response variable (with a potentially unbounded range). The canonical loss function for this problem (see Equation (\ref{rankparopt}) below) is not even continuous, though, so not immediately accessible for Gradient Boosting. \ \\

On the other hand, attracted by the immense computational efficiency of $L_2-$Boosting, it is tempting to somehow make $L_2-$Boosting accessible for model selection even to this non-standard loss. This idea also is suggested by the fact that $L_2-$Boosting (\cite{bu03}, \cite{bu06}) already approximates an optimal scoring rule for the hard continuous ranking problem, so why can't one just take the $L_2-$Boosting solution as solution for the ranking problem?  \ \\

Once we perform sparse model selection, structural risk minimization amounts to approximating the empirical risk minimizer with a model based on a sparse subset of the columns of the regressor matrix. More precisely, the final model often consists of a subset of the true relevant regressors which is determined by a sparse learning algorithm based on the respective loss function (possibly fruitfully enhanced by a Stability Selection, see \cite{bu10}, \cite{hofner15}). The fundamental question that arises is how we could guarantee that $L_2-$Boosting indeed selects all variables relevant for the ranking problem. \ \\

The answer is that we cannot. This statement is the core of this paper where we introduce a so-called ''column measure'' in order to exactly describe model selection and to provide a mathematical formulation of issues like that some sparse algorithm w.r.t. a loss function $L$ does not select all variables that were relevant for another loss function $\tilde L$. \ \\ 

Since the true column measure w.r.t. some loss function is not known, we will make key assumptions for the rest of this paper on its nature and of approximating properties of suitable model selection algorithms. We also identify resampling procedures that essentially assign weights to each row of a data matrix with the concept of a ''row measure''. The column measure, in combination with the row measure, defines a unified framework in which all learning procedures that invoke variable selection can be embedded (see \cite{TWphd} for a further discussion). \ \\

The rest of this paper is organized as follows. In Section \ref{prelim}, we recapitulate and introduce the mathematical framework of Gradient Boosting and especially $L_2-$Boosting. We restate the definition of the ranking loss functions that we need as well as the idea to use surrogate losses and the corresponding difficulties.  Section \ref{colmeasec} is devoted to the column measure and starts with a standard rejection sampling strategy that we already mentioned. We introduce the column measure and make several assumptions which we need throughout this work. Generalizing this idea, we also define the row measure. Moreover, we reformulate $k-$Step estimators in this framework and define a new estimator, the Expected $k-$Step estimator. Section \ref{singpartssec} provides the exact mathematical description of the issue that relevant variables which should be selected when concerning one loss function may happen to never be selected when using some other loss function. We show which consequences this problem has and which impact it causes on rejection sampling strategies. With all this conceptual infrastructure at hand we indeed render $L_2-$Boosting accessible for non-smooth losses, giving rise to a new algorithm SingBoost which somewhat pinpointedly could be termed a \textit{gradient-free functional Gradient Boosting algorithm}. We can show that under a Corr-min condition that ensures that each chosen variable is at least slightly correlated with the current residual, SingBoost enjoys estimation and prediction consistency properties.

\section{Preliminaries}\label{prelim} \

In this work, we always have data $\D=(X,Y) \in \R^{n \times (p+1)}$ where $Y_i \in \mathcal{Y} \subset \R$ and $X_i \in \mathcal{X} \subset \R^p$ where $X_i$ denotes the $i-$th row and $X_{\cdot,j}$ the $j-$th column of the regressor matrix $X$. \ \\ \ \\

\subsection{Boosting} \

The following generic functional Gradient Boosting algorithm goes back to Friedman (\cite{friedman01}). We use the notation of \cite{bu07}.\index{Boosting!Gradient Boosting|(} \ \\  \ \\

\begin{algorithm}[H] 
\label{funcBoosting}
\textbf{Initialization:} Data $(X,Y)$, step size $\kappa \in ]0,1]$, number $m_{iter}$ of iterations and \begin{center} $ \displaystyle \hat f^{(0)}(\cdot)\equiv \argmin_c\left(\frac{1}{n}\sum_i L(Y_i,c)\right) $ \end{center} as offset value\;
\For{$k=1,...,m_{iter}$}{
Compute the negative gradients and evaluate them at the current model: \begin{center} $ \displaystyle U_i=-\partial_f L(Y_i,f)|_{f=\hat f^{(k-1)}(X_i)} $ \end{center} for all $i=1,...,n$\;
Treat the vector $U=(U_i)_i$ as response and fit a model \begin{center} $ \displaystyle (X_i,U_i)_i \overset{\text{base procedure}}{\longrightarrow} \hat g^{(k)}(\cdot) $ \end{center} with a preselected real-valued base procedure\;
Update the current model via \begin{center} $ \displaystyle \hat f^{(k)}(\cdot)=\hat f^{(k-1)}(\cdot)+\kappa \hat g^{(k)}(\cdot) $  \end{center} 
}
\caption{Generic functional Gradient Boosting} 
\end{algorithm} \ \\ 

The function $L$ is considered to be a loss function with two arguments. One is the response, the other one is the predicted response, understood in a functional way as the model that is used for prediction. The loss function has to be differentiable and convex in the second argument. The main idea behind this algorithm is to iteratively proceed along the steepest gradient.  \ \\

As pointed out in \cite{bu07}, the models $\hat g^{(k)}$ can be regarded as an approximation of the current negative gradient vector.\index{Boosting!Gradient Boosting|)} \ \\

In this work, we concentrate on Gradient Boosting for linear regression with square loss, that is \textbf{Least Squares Boosting}, denoted by $L_2-$Boosting here. This variant has been shown to be estimation and prediction consistent even for very high dimensions (\cite[Thm. 1]{bu06}, \cite[Thm. 12.2]{bu}). It can be computed in an extremely efficient way using a componentwise linear regression procedure as described in Algorithm \ref{l2Boosting} (cf. \cite{bu06}). \ \\

The update step has to be understood in the sense that every $j-$th component of $\hat \beta^{(k)}$ remains unchanged for $j \notin \{1, \hat j_k\}$ and that only the intercept and the $(\hat j_k)-$th component are modified (if an intercept is fitted, we have $X_{\cdot, 1}=1_n$). Note that the formula to compute the residuals indicates that one is interested in following the steepest gradient. This is true in componentwise $L_2-$Boosting since for the standardized loss function \vspace{0.1cm} \begin{center} $ \displaystyle L_2(y,f):=\frac{1}{2}(y-f)^2, $ \end{center} the negative gradient w.r.t. $f$ is given by $(y-f)$.  \ \\ \ \\

\begin{algorithm}[H] 
\label{l2Boosting}
\textbf{Initialization:} Data $(X,Y)$, step size $\kappa \in ]0,1]$, number $m_{iter}$ of iterations and parameter vector $\hat \beta^{(0)}=0_{p+1}$\;
Compute the offset $\bar Y$ and the residuals $r^{(0)}:=Y-\bar Y$\;
\For{$k=1,...,m_{iter}$}{
\For{$j=1,...,p$}{
Fit the current residual $r^{(k-1)}$ by a simple least squares regression model using the predictor variable $j$\;
Compute the residual sum of squares\;
}
Take the variable $\hat j_k$ whose simple model $\hat \beta_{\hat j_k} \in \R^{p+1}$ provides the smallest residual sum of squares\;
Update the model via $\hat \beta^{(k)}=\hat \beta^{(k-1)}+\kappa \hat \beta_{\hat j_k}$\;
Compute the current residuals $r^{(k)}=Y-X\hat \beta^{(k)}$
}
\caption{Componentwise least squares or $L_2-$Boosting} 
\end{algorithm} \ \\ \ \\

Note that $L_2-$Boosting can be considered as a (generalized) rejection sampling strategy. In each step, one proposed $p$ different simple linear regression models based on the loss function $L_2$ and only accepts the model whose training $L_2-$loss is minimal. \ \\

\begin{rem} In practice, one computes the correlation of each column with the current residual since the linear regression model which improves the current combined model most, i.e., whose training $L_2-$loss is minimal, is always the linear regression model based on the column which has the highest absolute correlation with the current residual. This has been shown in a previous version of \cite{zhao07}, i.e., \cite{zhao04}. This is preferable as correlations can be updated very efficiently in this case. \end{rem} \ \\

\subsection{Empirical risk minimization and ranking} \

This subsection compiles notions and notation to cover ranking problems in the context of statistical learning. More specifically, one distinguishes several types of ranking problems. First of all, one can differentiate between \textbf{hard} (ordering of all instances), \textbf{weak} (identifying the top $K$ instances) and \textbf{localized} ranking problems (like the weak one, but also returning the correct ordering), see \cite{clem08b}. We head for hard ranking. \ \\

The second dimension concerns the scale of the response variable. If $Y$ is binary-valued, w.l.o.g. $\mathcal{Y}=\{-1,1\}$, then a ranking problem that intends to retrieve the correct ordering of the probabilities of the instances to belong to class 1 is called a \textbf{binary or bipartite ranking problem}\index{Ranking problems!Bipartite|(}. If $Y$ can take $d$ different values, a corresponding ranking problem is referred to as a \textbf{$d-$partite ranking problem}\index{Ranking problems!$d-$partite|(} and for continuously-valued responses, one speaks of a \textbf{continuous ranking problem}\index{Ranking problems!Continuous|(}.  \ \\

While there exist a vast variety of machine learning algorithms to solve binary ranking problems (e.g., \cite{clem08c}, \cite{clem08d}, \cite{joachims02}, \cite{herb}, \cite{pahi}, \cite{freund}), the only approach for continuous ranking problems that we are aware of is a tree-type approach in \cite{clem18}. \cite{freund} invoked a convex surrogate loss with a pairwise structure to develop the RankBoost algorithm which solves the hard binary ranking problem, but a similar strategy is not meaningful for continuous ranking problems, as we have discussed in a review paper (\cite{TW19b}). \ \\

Returning to our motivating example from tax auditing, it turns out that it is helpful to account for the severity of the damage and hence one is led to the continuous ranking problem. It even pays off to model severity with an unbounded range.  \ \\

For this purpose, empirical risk minimization requires the definition of a suitable risk function. We head for the hard ranking problem. In this context, we take over the loss functions as introduced in \cite{clem05}. One starts with the class $\mathcal{R}$ of ranking rules $r: \mathcal{X} \times \mathcal{X} \rightarrow \{-1,1\}$ where $r(X,X')=1$ indicates that $X$ is ranked higher than $X'$. Then, for ranking rule $r \in \mathcal{R}$, the corresponding risk is given by \vspace{0.25cm} \begin{center} $ \displaystyle  R^{hard}(r):=\E[I((Y-Y')r(X,X')<0)] $ \end{center} \vspace{0.25cm} where $X'$ resp. $Y'$ is an identical copy of $X$ resp. $Y$. In fact, this is nothing but the probability of a misranking of $X$ and $X'$. Thus, empirical risk minimization intends to find an optimal ranking rule by solving the optimization problem \vspace{0.25cm} \begin{center} $ \displaystyle \min_{r \in \mathcal{R}}\left(L_n^{hard}(r)=\frac{1}{n(n-1)}\mathop{\sum \sum}_{i \ne j} I((Y_i-Y_j)r(X_i,X_j)<0)\right). $ \end{center} \vspace{0.25cm} For the sake of notation, the additional arguments in the loss functions are suppressed. Note that $L_n^{hard}$, i.e., the hard empirical risk, is also the hard ranking loss\index{Ranking losses!Hard|(} function which reflects the global nature of hard ranking problems.  In practice, one defines some scoring function $s_{\theta}: \mathcal{X} \rightarrow \R$ for $\theta \in \Theta \subset \R^p$. Then, the problem carries over to the parametric optimization problem \vspace{0.25cm} \begin{equation} \label{rankparopt} \min_{\theta \in \Theta}\left(L_n^{hard}(\theta)=\frac{1}{n(n-1)}\mathop{\sum \sum}_{i \ne j} I((Y_i-Y_j)(s_{\theta}(X_i)-s_{\theta}(X_j))<0)\right). \end{equation}\vspace{0.25cm}Let us emphasize that the hard ranking loss function is far away from being a suitable candidate loss function for Gradient Boosting. Not only does it fail to be differentiable, it is not even continuous and due to the non-convexity, even subgradients do not exist. \ \\

\begin{rem} As discussed in more detail in \cite{TW19b}, the straight-forward approach to use surrogate losses (as successfully employed in RankBoost for bipartite ranking, cf. \cite{freund}) is \textbf{not} a remedy for continuous ranking problems in general, in particular, if the range of the response is unbouded. The respective strategies may lead to nearly empty model selections where in fact good predictors exist, it is likely to cause numerical problems and so forth. \end{rem} 

\section{Column measure} \label{colmeasec} \

\subsection{Na\"{i}ve strategy} \
    
According to \cite[Prop. 1]{clem18}, the conditional expectation $\E[Y|X]$ is an optimal scoring rule for the hard continuous ranking problem. Therefore, using the fact that $L_2-$Boosting computes a sparse approximation of $\E[Y|X]$, a simple strategy could be to use the $L_2-$Boosting model as solution of the ranking problem. Somewhat more sophisticated, one would draw subsamples from the data and to evaluate the out-of-sample performance w.r.t. the hard ranking loss of the model based on each subsample. This could be seen as a (generalized) rejection sampling strategy since only the best model was taken. \ \\

This idea leads to the following fundamental question: Why should the optimal sparse regression model, i.e., w.r.t. the squared loss, be identical to the optimal sparse hard ranking model w.r.t. the loss function $\tilde L$ in Equation (\ref{rankparopt}) according to the set of selected variables? \ \\

The definition of the column measure that we introduce in this section provides an exact mathematical description of the problem that the optimal sparse models may differ. \ \\

\subsection{The definition of a column measure} \

\begin{Def} \label{colmeasure} Suppose we have a data set $\D:=(X,Y) \in \R^{n \times (p+1)}$ and a loss function $L: (\mathcal{Y} \times \mathcal{Y}) \rightarrow \R_{\ge 0}$ that should be minimized empirically. Then the \textbf{column measure w.r.t. $L$}\index{Column measure} is defined as the map \begin{center} $ \displaystyle \nu^{(L)}: (\{1,...,p\}, \mathcal{P}(\{1,...,p\})) \rightarrow ([0,p], \mathbb{B} \cap [0,p]), $ \end{center}\begin{center} $ \displaystyle  \nu^{(L)}: \{j\} \mapsto \nu^{(L)}(\{j\})=:\nu_j^{(L)} \in [0,1] \ \forall j \in \{1,...,p\} $ \end{center} and assigns an importance to each predictor. \end{Def}

\begin{rem} Note that the column measure is not a probability measure. However, for each singleton $\{j\} \in \{1,...,p\}$, the quantity may be seen as a probability for the corresponding predictor to be chosen. Of course, in the case that $p \rightarrow \infty$, we do no longer have a finite measure but still a $\sigma-$finite measure. \end{rem}

\begin{rem} We used this definition primarily because of the easy interpretation\index{Column measure!Interpretation} that the value assigned to each singleton can be interpreted as a plausibility with values in $[0,1]$. For theoretical purposes, it may be convenient to normalize the column measure to get a probability measure for finite $p$. Thus, one would just divide each $\nu^{(L)}(\{j\})$ by $\nu^{(L)}(\{1,...,p\})$. \end{rem} \

Of course, we do not know the true column measure $\nu$ for any loss function. But we will make the following assumptions throughout this paper. \ 

\begin{ass}[\textbf{Available predictors for varying sample sizes}] \label{availass} We cover two situations, i.e., a) a fixed set of available predictors and b) a set of available predictors varying in $n$. For the respective situation we assume 

\textbf{a)} We work with a fixed predictor set $\mathcal{X}$ in the sense that the column set does not change in $n$, only new rows get attached, and the number of available predictors $p_n\equiv p$. 

\textbf{b)} The available set of predictors may grow in $n$, but for each predictor index $j$ there is a sample size $N=N(j)$ such that this predictor is available for each sample size $n>N$. \end{ass} 

\begin{ass}[\textbf{Weak convergence of column measures}] \label{nuass} Let $L$ be a loss function as before and let $\hat \nu_n^{(L)}$ (with $\hat \nu_{n,j}^{(L)}:=\hat \nu_n^{(L)}(\{j\})$) be the empirical column measure\index{Column measure!Empirical} that has been computed by a variable selection consistent, $L-$adapted learning procedure. Then we assume that \begin{center} $ \displaystyle \hat \nu_n^{(L)} \overset{w}{\longrightarrow} \nu^{(L)} $. \end{center} \end{ass}

\begin{rem}[\textbf{Weak convergence}]\index{Column measure!Weak convergence} Note that by the Portmanteau theorem (see e.g. \cite[Thm. 4.10]{elstrodt}), weak convergence stated in Assumption \ref{nuass} is equivalent to the condition \begin{center} $ \displaystyle \hat \nu(B) \longrightarrow \nu(B) $ \end{center} for all $\nu-$continuity sets $B$ which in situations a) and b) of our Assumption \ref{availass} is just pointwise convergence \begin{center} $ \displaystyle \hat \nu^{(L)}_{n,j} \overset{n \rightarrow \infty}{\longrightarrow} \nu_j^{(L)} \ \ \forall j \in \{1,...,p\}. $ \end{center} \end{rem}  

\begin{rem} Assumption \ref{nuass} makes no statement as to variable selection consistency\index{Column measure!and variable selection consistency} since this property only means that the true model is chosen asymptotically. This has to be understood in the language of column measures that asymptotically, we have \begin{equation} \begin{gathered} \hat \nu_{n,j}^{(L)} \longrightarrow 0 \ \forall j \notin S^0 \\ \hat \nu_{n,j}^{(L)} \longrightarrow c_j^{(L)}>0 \ \forall j \in S^0 \end{gathered} \end{equation} for $n \rightarrow \infty$ where $S^0$ is the true predictor set. \end{rem} 

\begin{rem} Working with weak convergence from Assumption \ref{nuass}, we potentially could extend the definition of column measures on subsets of $\N$ to ''column measures'' on subsets of $\R$ or other uncountable sets. For instance, assume that we had time-dependent measurements where the time index $t$ can be treated as element of the uncountable interval $[0,T]$ for some $T \in \R_{>0}$ and the goal is to represent this sequence of measurements by a sparser sequence. Then we essentially had to work with ''column measures'' on $\R$. \end{rem}  
 \ \\

\subsection{Randomness and Column Measure} \

From an accurate statistical modelling point of view, exactly assessing the random mechanism behind the observations, the column measure could be counter-intuitive, because model selection is a matter of design and not of randomness.  \ \\

From a mathematical point of view though, the column measure formalizes an aggregation mechanism exactly parallel to the usual construction of a measure -- in this case on the discrete sample space of the column indices.      \ \\

In particular, all concepts from measure theory carry over and lead to a both fruitful and powerful language and theoretical infrastructure, which we want to make extensive use of. This approach parallels the introduction of martingale measures in Mathematical Finance, where conceptually one starts from non-random sets of cash flows which through aggregation are then endowed with corresponding probability concepts, compare the exposition in \cite{irle}. The most striking advantages of the present column measure approach, hence of this language and infrastructure are, in our opinion \begin{itemize}
\item a conceptually uniform and consistent treatment of different column selection mechanisms
\item  interpretability as particular rejection sampling schemes (with the theory of theses schemes behind it)
\item  interpretability of our algorithm SingBoost (see Algorithm \ref{singboost}) as a particular change-of-measure operation (with all underlying theory for the latter)
\item  ready availability of convergence concepts (with limit theorems like Law of Large Numbers and Central Limit Theorems)
\end{itemize}
\ \\

\subsection{The row measure} \

From a computer science point of view, the empirical starting point in any learning algorithm is a data matrix (or, more accurately, a data frame in R language, allowing for varying data types in the columns). As any matrix has two indices for rows and columns, it is evident that our column measure by trivial duality can be translated to a row measure, i.e., a measure acting on the observations, just as well. In fact, the corresponding row measure is nothing but the well known and well studied empirical measure. Still, in computational algorithms for aggregation, some added value in interpretation is obtained by formulating aggregation steps through the row-measure perspective. We spell this out here. \ \\

\begin{Def} \label{rowmeasure}\index{Row measure} Suppose we have $n$ observations $X_1,...,X_n \in \mathcal{X}$ where $\mathcal{X}$ is some space. W.l.o.g., we think of $X:=(X_1,...,X_n)$ as an $(n \times 1)-$dimensional column vector or, in the case of multivariate observations, as matrix with $n$ rows. Then we define the \textbf{row measure} as the map \begin{center} $ \displaystyle \zeta: (\{1,...,n\}, \mathcal{P}(\{1,...,n\})) \rightarrow ([0,1], \B \cap [0,1]), $ \end{center} \begin{center} $ \displaystyle \zeta: \{i\} \mapsto \zeta(\{i\})=:\zeta_i \in [0,1] \ \forall i \in \{1,...,n\} $ \end{center} which assigns a weight to each row. \end{Def} \ 

In other words, special sampling techniques replace the uniform distribution on the row indices, i.e., the \textbf{uniform row measure} to some other measure which is tailored to those data sets. \ \\

As we have seen, \textbf{empirical column measures that are able to assign values to singletons in the whole interval $[0,1]$ can only be computed by aggregation of different models}. If these different models were fitted on different subsamples or Bootstrap samples of the data, then we actually induced an empirical column measure in the sense of the following definition. \ \\

\begin{Def} \label{indcolmea} Let $\D:=(X,Y) \in \R^{n \times (p+1)}$ be a data matrix and let $\zeta$ be a row measure on the row indices. Assume that a fitting procedure that optimizes some empirical loss corresponding to a loss function $L: \mathcal{Y} \times \mathcal{Y} \rightarrow \R_{\ge 0}$ is given. Let $B$ be the number of samples drawn from $\zeta$. Let $(\hat \nu^{(L)})^{(b)}$ be the empirical column measure fitted on the data set reduced to the rows given by the $b-$th realization of $\zeta$. Then we call the resulting empirical column measure $\hat \nu^{(L)}(\zeta)$ with \begin{center} $ \displaystyle \hat \nu^{(L)}(\zeta)(\{j\}):=\frac{1}{B}\sum_b (\hat \nu^{(L)})^{(b)}(\{j\}) $ \end{center} the \textbf{empirical column measure induced by the row measure $\zeta$}\index{Column measure!Induced by a row measure}. \end{Def}

\begin{rem} In fact, aggregating over all $B$ samples autonomizes against the concrete realizations of the row measure $\zeta$. With a little abuse of notation, one can see this \textbf{empirical column measure induced by $\zeta$ as empirical expectation of $\hat \nu^{(L)}$ w.r.t. $\zeta$}. \end{rem} 

\begin{rem} Without calling them measures, \cite{bu10} implicitly spoke of induced column measures when emphasizing that for a given penalty parameter $\lambda$, the set of selected predictors w.r.t. $\lambda$ also depends on the subsample. \end{rem}\

The obvious most general concept would then combine rows and columns, and work with a joint aggregation over both rows and columns, leading to a joint column- and row measure, assigning aggregation weights to individuall matrix cells. This is helpful as a concept to treat uniformly and consistently cell-wise missings and outliers as well as resampling schemes with cell individual weights (see e.g. \cite{alqallaf}, \cite{leung16}, \cite{rous16}). \ \\ \ \\

\subsection{Expected $k-$Step} \

The main difficulty up to this point is the question how empirical column measures can be incorporated into model fitting and prediction. As their corresponding empirical column measure is the ''$0/1-$limit case'', any existing model selection method requires at least the decision whether a predictor variable enters the final model or not. This is especially the case if one applies a Stability Selection (\cite{bu10}, \cite{hofner15}) where one indeed computes an empirical column measure by averaging several 0/1-valued column measures resulting from applying the respective algorithm, e.g. Boosting or a Lasso, on subsamples of the data. However, at the end, one transforms it again into a 0/1-valued empirical column measure using a pre-defined cutoff such that only the predictors with a sufficiently high empirical selection frequency are selected for the final model. \ \\

We try to explicitly account for the relative importance of each variable. Therefore, we present an estimator that is based on $k-$Step estimators on whose improvement property we rely. \ \\ 

We start with the important definition of influence functions that can be found in \cite{hampel} and was originally stated in \cite{hampel68}. \ \\ 

\begin{Def}\label{IC} Let $X$ be a normed function space containing distributions on a probability space $(\Omega, \mathcal{A})$ and let $\Theta$ be a normed real vector space. Let $T: X \rightarrow \Theta$ be a statistical functional. The \textbf{influence curve} or, more generally, \textbf{influence function} of $T$ at $x$ for a probability measure $P$ is defined as the derivative \begin{center} $ \displaystyle \IC(x,T,P):=\lim_{t \rightarrow 0}\left( \frac{T((1-t)P+t \delta_x)-T(P)}{t} \right)=\partial_t \left[T((1-t)P+t \delta_x) \right] \bigg|_{t=0} $ \end{center} where $\delta_x$ denotes the Dirac measure at $x$. \end{Def} \

This is just a special G\^{a}teaux derivative (see e.g. \cite{averbukh}) with $h:=\delta_x-P$. The influence curve can be regarded as an estimate for the infinitesimal influence of a single observation on the estimator. \ \\

In the setting of a smooth parametric model $\mathcal{P}=\{P_{\theta} \ | \ \theta \in \Theta\}$ for some $p-$dimensional parameter set $\Theta$, we assume that $\mathcal{P}$ is $L_2-$differentiable in some $\theta_0$ in the open interior of  $\Theta$ with $L_2-$derivative $\Lambda_{\theta_0}$ (\cite[Def. 2.3.6]{rieder}). Then the family $\Psi_2(\theta_0)$ of influence curves in $P_{\theta_0}$ is defined as the set of all open maps $\eta_{\theta_0}$ that satisfy the conditions \ \\

 \textbf{i)} $\eta_{\theta_0} \in L_2^p(P_{\theta_0})$, \ \ \ \textbf{ii)} $\E_{\theta_0}[\eta_{\theta_0}]=0$, \ \ \ \textbf{iii)} $\E_{\theta_0}[\eta_{\theta_0} \Lambda_{\theta_0}^T]=I_p $ \ \\  

where $I_p$ denotes the identity matrix of dimension $p \times p$. See \cite[Sec. 4.2]{rieder} for further details. \ \\

A general recipe to obtain estimators with a prescribed influence function is given by the one-step principle. I.e., to a $\sqrt{n}-$consistent starting estimator $\hat \theta_n$ one defines the One-step estimator $S_n^1$ by \begin{equation} \label{onestep} S_n^1:=\hat \theta_{n}+\frac{1}{n} \sum_{i=1}^n \eta_{\hat \theta_n}(x_i) \end{equation} is referred to as One-Step estimator. \ \\

In fact, from \cite[Thm. 6.4.8]{rieder}, $S_n^1$ indeed admits influence curve $\eta_{\theta}$: \

\begin{prop} \label{prop} In the parametric model $\mathcal{P}$, for influence function $\eta_{\theta}$ according to i)-iii), the One-step estimator defined in Equation (\ref{onestep}) is asymptotically linear in the sense that \begin{equation} S_n^1=\theta_0+\frac{1}{n} \sum_{i=1}^n \eta_{\hat \theta_0}(x_i)+R_n \end{equation} for some remainder $R_n$ such that $\sqrt{n}R_n \rightarrow 0$ stochastically. \end{prop}  \

See \cite[Sec. 6]{rieder} for further insights on this principle. The One-Step estimator can be identified with a first step in a Newton iteration as described in \cite{bickel75} for the linear model. In fact, this approximation amounts to usage of a functional Law of Large Numbers which in its weak/stochastic variant is an easy consequence of the routine functional Central Limit Theorem argument employed to derive asymptotic normality for M-estimators (cf. \cite[Cor. 1.4.5]{rieder}, \cite[Sec. 5.3+5.6]{vaart98}).  \ \\

Iterating this principle $k$ times with the estimate from the previous step as corresponding starting estimator, one obtains a $k-$Step estimator. \ \\

The adaption of $k-$Step estimation to model selection\index{$k-$Step estimator} is straightforward. For notational simplicity, in the sequel of this subsection, we stick to Assumption \ref{availass} a) but note that up to a more burdensome notation we could equally cover the case of Assumption \ref{availass} b). Having selected a set $J \subset \{1,...,p\}$ of columns, w.l.o.g. ordered in an ascending sense, the starting estimator is a $\sqrt{n}-$consistent estimator which is essentially based on the reduced data set $(X_{\cdot,J},Y)$. This is the usual procedure for estimation. We just have to account for the influence curve. \ \\

In fact, by the previously performed model selection, we just have to estimate the reduced parameter $\theta_J$. This reduction can be thought of the mapping \begin{center} $ \displaystyle \tau: \R^p \rightarrow \R^{|J|}, \ \ \ \theta \mapsto (\theta_j)_{j \in J}. $ \end{center} Since $\tau$ is differentiable and $|J|=:q \le p$, we can compute the partial influence curve \begin{center} $ \displaystyle \partial_{\theta} \tau(\theta)\eta_{\theta}=(\eta_{\theta})_J. $ \end{center} which is already admissible if $\eta_{\theta} \in \Psi_2(P_{\theta})$ due to \cite[Def. 4.2.10]{rieder}. 

Thus, the One-Step estimator for $\theta_J$ is given by \begin{equation} \label{onestepred} S_{n,J}^1=\widehat{(\theta_J)_n}+\frac{1}{n}\sum_i (\eta_{\widehat{(\theta_J)_n}})_J(X_i,Y_i).  \end{equation} The extension to the $k-$Step estimator is straightforward. Formally, given the $k-$Step estimator $S_{n,J}^k$ for the parameter $\theta_J$, a resulting $p-$dimensional parameter is clearly given by any preimage \begin{center} $ \displaystyle S_n^k=\tau^{-1}(S_{n,J}^k). $ \end{center}    \ Note that during the one-($k-$)step estimation as in Equation (\ref{onestepred}), the set of selected variables $J$ remains unchanged. \ 

\begin{rem}\label{picrem} Assume for simplicity that $J=\{1,...,q\}$. Then, the so-called partial influence curve $(\eta_{\theta})_J$ satisfies \begin{center} $ \displaystyle \E[(\eta_{\theta})_J \Lambda_{\theta}^T]=(I_q,0_{p-q}) \in \R^{q \times p} $ \end{center} which indeed is obviously true since the last $(p-q)$ rows of $I_q$ just have been dropped. See \cite{rieder} for further details on partial influence curves. \end{rem}  \

Observe the following fact: We have done model selection and have gotten a subset $J \subset \{1,...,p\}$ of relevant variables. In the ''$0/1-$setting'', this corresponds to the column measure $\hat \nu^{(L)}=(I(j \in J))_j$. Due to Equation (\ref{onestepred}), the One-Step $S_n^1$ is asymptotically equal to \vspace{0.25cm} \begin{equation} \label{expkstep} \E_{\hat \nu^{(L)}}\left[\hat \theta_n+\frac{1}{n}\sum_i \eta_{\hat \theta_n}(X_i,Y_i)\right]:=\left(\left((\hat \theta_n)_j+\frac{1}{n}\sum_i (\eta_{\hat \theta_n}(X_i,Y_i))_j\right)I(j \in J)\right)_{j=1}^p, \end{equation} where $\hat \theta_n$ is the estimated coefficient on the complete data set. This is true since the starting estimator is $\sqrt{n}-$consistent, so the difference in the components of the parameter corresponding to $J$ gets negligible. The same is true for the influence curve by theorem \cite[Thm. 6.4.8a)]{rieder}. The components corresponding to $J^c$ are zero in both cases. Implicitly, the One-Step estimator that is computed after model selection is asymptotically an expectation of the ''original'' One-Step estimator (i.e., without model selection) w.r.t. a suitable column measure taking values in $\{0,1\}$. Note again that this leads to partial influence curves. \

\begin{Def} \label{kstepdef} Let $\D:=(X,Y) \in \R^{n \times (p+1)}$ be a data matrix. Let $\hat \nu_n$ be some empirical column measure on the set $\{1,...,p\}$ and let $S_n^1$ be a One-Step estimator based on $n$ observations. Then we define the \textbf{Expected One-Step estimator}\index{$k-$Step estimator!Expected} \begin{equation} \label{exponestep} \E_{\hat \nu_n}[S_n^1]:=\hat \nu_{n,vec} \circ_H S_n^1 \end{equation} where $\circ_H$ is the component-wise (Hadamard) product and \begin{center} $ \displaystyle \hat \nu_{n,vec}:=(\hat \nu_n(\{j\}))_{j=1}^p. $ \end{center} \end{Def}

\begin{rem} The Expected One-Step estimator can be thought of weighting each component of the coefficient by its empirically estimated relevance. Since these weights only take values in $[0,1]$, we essentially perform a \textbf{shrinkage} of the coefficients that are not chosen with an empirical probability of one. \end{rem}

\begin{rem} The usual One-Step estimator is constructed to get an estimator with a given influence curve. The Expected One-Step estimator can be seen as an extension in the sense that the resulting estimator both has the desired (partial) influence curve and is based on a given column measure. \end{rem} \

\subsection{Application in resampling} \

The concept of the Expected One-($k-$)Step estimator as introduced in Definition \ref{kstepdef} turns out particularly useful in connection with resampling schemes. 

Due to the fact that stochastic convergence is not convex in the sense that for a triangular scheme $(X_{i;n})_{i=1,\ldots,n;\,n\in\N}$, it may happen that $\mean_{i=1,\ldots,n} X_{i;n} \nrightarrow 0$ for $n\rightarrow\infty$  while still $X_{i;n}\rightarrow 0$ stochastically for $n\rightarrow\infty$ and for each $i$, we have to sharpen the remainder $R_n$ in Proposition~\ref{prop} such that \begin{equation} \label{sharpass} \| \,\sqrt{n}\,R_n\|_{L_1(P_{\theta_0})}\longrightarrow 0. \end{equation} In fact, convergence (\ref{sharpass}) is a matter of uniform integrability and can be warranted in many cases by suitable exponential concentration bounds (see \cite[Prop. 2.1b)+Thm. 4.1]{ruck10}).   \ \\

We now investigate the problem if we can combine the concept of $k-$Step estimators with any column measure. Assume that we already have drawn $B$ samples from the data and computed the One-Step estimator $\hat \theta^{(b)}$ for each $b=1,...,B$. Then averaging leads to the estimator \begin{equation} \label{avgonestep} \hat \theta^{(avg)}:=\frac{1}{B}\sum_b \hat \theta^{(b)}+\frac{1}{n}\sum_i \frac{1}{B} \sum_b \eta_{\hat \theta^{(b)}}(X_i,Y_i) \end{equation} which by (\ref{sharpass}) and Proposition \ref{prop} still is asymptotically linear with influence function $\eta_{\theta_0}$.  What does actually happen? Assume there is a subset $K \subset \{1,...,B\}$ such that a variable $j_0$ is chosen in all iterations $b \in K$. Then the $j_0-$th component of the estimated coefficient in Equation (\ref{avgonestep}) would be \begin{equation} \label{avgonestepcomp} \hat \theta^{(avg)}_{j_0}=\frac{1}{B}\sum_{b \in K} \hat \theta_{j_0}^{(b)}+\frac{1}{n}\sum_{i=1}^n \sum_{b \in K}(\eta_{\hat \theta^{(b)}}(X_i,Y_i))_{j_0}. \end{equation} Note that the empirical column measure is $\hat \nu(\{j_0\})=\frac{|K|}{n}$ in this case. \ 

\begin{thm} Assume that the variable selection procedure that leads to the empirical column measure $\hat \nu_n$ is variable selection consistent and that the learning procedure that leads to the estimated coefficients is $\sqrt{n}-$consistent. Then the Expected One-Step estimator is variable selection consistent\index{Variable selection consistency!of the expected One-Step} and is $\sqrt{n}-$estimation consistent for the coefficients $j$ with $\hat \nu_n(\{j\})=1 \ \forall n$ assuming (\ref{sharpass}). \end{thm}

\begin{bew} \textbf{a)} By assumption, the variable selection procedure is variable selection consistent (see e.g. \cite{bu}), so $\hat S=S^0$ asymptotically. Then the empirical column measure is, asymptotically, zero in all entries $j \notin S^0$, so the Expected One-Step estimator is already correct in these components. The learning procedure is $\sqrt{n}-$consistent, so any estimated coefficient $\hat \theta^{(b)}$ falls into an $o_P(n^{-1/2})-$neigh\-borhood of the true coefficient, hence the arithmetic mean over all $b$ also does. That implies, as for the initial estimator, it does not depend if it is chosen as arithmetic mean over all $\hat \theta^{(b)}$ or if it is computed once since its distance to the true coefficient gets negligible for $n \rightarrow \infty$. The same holds for the arithmetic mean of influence curves due to property i) in theorem \cite[Thm. 6.4.8a)]{rieder}, provided that the variable is chosen with an empirical probability of one as assumed. 

\textbf{b)} By the variable selection consistency, a growing number of observations guards against false negatives. Thus, asymptotically this case does not have to be treated. The same is true for false positives.  

\begin{flushright} $_\Box $ \end{flushright} \end{bew}

\section{''Likelihood'' ratio of two column measures} \label{singpartssec} \

\subsection{Singular parts} \

When performing Gradient Boosting w.r.t. different loss functions, it becomes evident that the selection frequencies of the variables differ from loss function to loss function. Therefore, we cannot assume that the (empirical) column measure $\hat \nu^{(L)} $ is identical to the (empirical) column measure $\hat \nu^{(\tilde L)}$ for two different loss functions $L$ and $\tilde L$. So why should equality be even true for the true column measures $\nu^{(L)}$ and $\nu^{(\tilde L)}$? Ideas related to the rejection sampling strategy will extend to the situation where these two measures agree in terms of which variables are important and only differ in the concrete selection probabilities; mathematically this is nothing but the statement that these two measures are equivalent in the sense of domination or, for column measures depending on sample size $n$, of contiguity (cf. \cite{lecam86}). \ \\

An issue arises if Gradient Boosting is applied w.r.t loss function $L$ and evaluated the performance in the loss function $\tilde L$, and if some relevant variables for loss function $\tilde L$ are never selected by Gradient Boosting w.r.t. $L$. In this case, rejection sampling can never reliably succeed in finding all relevant variables for loss function $\tilde L$. This leads to the question if it is realistic to assume that suitable model selection algorithms for different loss functions select different variables.
\ \\

Let $S^0$ be the true set of relevant variables. Then an algorithm with the screening property (see \cite{bu}) tends to select a superset of $S^0$, so is prone to overfitting. If we compared two algorithms for different loss functions that both have the screening property, it asymptotically may happen that they select \textbf{different noise variables}. This would not be an issue as long as $S^0$ is contained in the selected sets of variables for both algorithms. If even the variable selection property is valid, then the true set $S^0$ is chosen asymptotically.  \ \\

So far, however, there is, to the best of our knowledge, no such asymptotic dominance statement of $\nu^{(L)}$ w.r.t. $\nu^{(\tilde L)}$ available. To be on the safe side, it does make sense to consider singular parts of $\nu^{(\tilde L)}$ w.r.t. $\nu^{(L)}$.  \



\begin{ass} \label{diffcolass} We assume that even the true column measures for different loss functions $L$ and $\tilde L$ differ, i.e., \begin{center} $ \displaystyle \nu^{(L)} \ne \nu^{(\tilde L)}. $ \end{center} \end{ass} \

This motivates to adapt the definition of domination and singular parts for our column measure framework. \

\begin{Def} \label{domcolmea} Let $L$, $\tilde L$ be loss functions as before and let $\D:=(X,Y) \in \R^{n \times (p+1)}$ be a data set. 

\textbf{a)} We call the column measure $\nu^{(L)}$ \textbf{dominated by}\index{Column measure!Domination} $\nu^{(\tilde L)}$, written $\nu^{(\tilde L)}\gg\nu^{(L)}$ if \begin{center} $ \displaystyle \nu^{(\tilde L)}(\{j\})=0 \Longrightarrow \nu^{(L)}(\{j\})=0 \ \forall j \in \{1,..,p\}. $ \end{center} 

\textbf{b)} We say that $\nu^{(L)}$ and $\nu^{(\tilde L)}$ are \textbf{equivalent}\index{Column measure!Equivalence}, $\nu^{(L)} \sim \nu^{(\tilde L)}$, if \begin{center} $ \displaystyle \nu^{(L)}\gg \nu^{(\tilde L)} \ \text{and} \ \nu^{(\tilde L)}\gg\nu^{(L)}, $ \end{center}  or just written differently, if \begin{center} $ \displaystyle \nu^{(L)}(\{j\})=0 \Longleftrightarrow \nu^{(\tilde L)}(\{j\})=0. $ \end{center}     

\textbf{c)} If $\nu^{(\tilde L)}\gg\nu^{(L)}$ and if the measures are not equivalent, then we call the greatest set $J_{\tilde L}^L \subset \{1,...,p\}$ such that \begin{center} $ \displaystyle  \nu^{(L)}(\{j\})=0 \ne \nu^{(\tilde L)}(\{j\}) \ \forall j \in J_{\tilde L}^L $ \end{center} the \textbf{singular part of $\nu^{(\tilde L)}$\index{Column measure!Singular parts} w.r.t. $\nu^{(L)}$}.  \ \\

\textbf{d)} In the situation  with column measures $\nu_n^{(\tilde L)}$ and $\nu_n^{(L)}$ depending on the sample size $n$,  replacing equality with $0$ by stochastic convergence to $0$ in a) and b), we obtain the notion of that $(\nu_n^{(L)})$ be \textbf{continguous} (\textbf{asymptotically equivalent}) to $(\nu_n^{(\tilde L)})$ in the sense of \cite{lecam86} which is denoted by the same symbols $(\nu_n^{(L)})\gg (\nu_n^{(\tilde L)})$ ($ (\nu_n^{(L)})\sim(\nu_n^{(\tilde L)})$). \end{Def} 

\begin{rem} One could extend the definition of the column measure by explicating the learning procedure that is used to perform risk minimization. \\
It certainly will make a difference if the quadratic loss is minimized by ordinary least squares or by $L_2-$Boosting where the latter will generally result in a much sparser model, generating a singular part, too. We do not think that this extension would be meaningful. If one is interested in the sparsity of the models computed by different algorithms for the same loss function, one does not really need the column measure. On the other side, of course the empirical column measure from the ordinary least squares model would ''dominate'' every other model's empirical column measure, but if every variable is selected anyways, the notion of column measure is not needed,  and in this case would lead trivially to assigning $1$ to each singleton. \end{rem}

\begin{ass} If we compare the relevant sets of predictors for different loss functions, we always assume that a suitable learning algorithm performing (sparse) variable selection has been applied without accounting for it specifically. \end{ass} 

Since we encounter discrete measures on a finite measurable space, we can easily find the Lebesgue decomposition of either measure if none of them dominates the other one. \

\begin{ex}[\textbf{Lebesgue decomposition}]\index{Column measure!Lebesgue decomposition} Assume we have a set $J$ where both measures are nonzero for every singleton contained in $J$. Additionally, we have singular parts $J_{\tilde L}^L$, $J_L^{\tilde L}$, so the union of all three sets is $\{1,...,p\}$. Then we trivially get the Lebesgue decompositions \begin{center} $ \displaystyle \nu^{(\tilde L)}=\nu^{(\tilde L)}|_{J_L^{\tilde L} \cup J}+\nu^{(\tilde L)}|_{J^L_{\tilde L}}, $ \end{center} \begin{center} $ \displaystyle \nu^L=\nu^L|_{J_{\tilde L}^L \cup J}+\nu^{(L)}|_{J^{\tilde L}_L}. $ \end{center}\ \\ Then $\nu^{(\tilde L)}|_{J^L_{\tilde L}} \perp \nu^{(L)}$ and $\nu^L|_{J^{\tilde L}_L} \perp \nu^{(\tilde L)}$. \end{ex}

\begin{ass} We even go further than Assumption \ref{diffcolass} and assume that the column measures $\nu^{(L)}$ and $\nu^{(\tilde L)}$ already differ on the set $S^0$ while not necessarily having singular parts, so we assume \begin{center} $ \displaystyle \nu^{(L)}|_{S^0} \ne \nu^{(\tilde L)}|_{S^0}, \ \ \ \text{i.e.,} \ \ \ \exists j \in S^0: \nu_j^{(L)} \ne \nu_j^{(\tilde L)}. $ \end{center}  \end{ass}  \

\subsection{Consequences to model selection}  \

With the definition of the column measure w.r.t. different loss functions and their potential singular parts, we once more take a look at the rejection sampling idea. 

Each proposed model (e.g. a Boosting model or a model determined by Stability Selection with Boosting or with the Lasso) provided a ''raw'' empirical column measure which is not necessarily stable, though. However, this measure corresponds to the selected loss function, i.e., the loss $L$. Performing rejection sampling by evaluating the model w.r.t. the loss function $\tilde L$ effectively picks one of these proposed raw column measures. The resulting column measure thus differs from the one that we would have gotten for the original target loss $\tilde L$. Mathematically spoken, through rejection sampling, we implicitly made a \textbf{change of measure}.  \ \\ 

In the same way as a change of measure requires that the initial measure dominates the target measure, we get a similar requirement for a rejection sampling strategy. 

\begin{rem}[\textbf{Rejection sampling and domination}]\index{Rejection sampling!Domination} A rejection sampling strategy where models are fitted w.r.t. a loss function $L$ and evaluated in another loss function $\tilde L$ is inadequate if \begin{equation} \label{rejsingpart} J_{\tilde L}^L \ne \emptyset, \end{equation} or in other words, if the measure $\nu^{(\tilde L)}$ has a singular part w.r.t. $\nu^{(L)}$.\ \\

In turn, our strategy can only be meaningful if \begin{equation} \label{rejdompart} \nu^{(L)}\gg\nu^{(\tilde L)}. \end{equation} Evidently, singular parts of $\nu^{(L)}$ w.r.t. $\nu^{(\tilde L)}$ can lead to false positives (w.r.t. $\tilde L$), but the rejection sampling strategy can still work if these false positives can be appropriately neglected in subsequent stabilization steps. \ \\

Considering the original rejection sampling (e.g. \cite[Alg. 3.4]{ripley}), mere dominance would not be sufficient, but instead one requires the likelihood ratios as bounded. In our setting of finitely many columns this is no extra requirement though, as boundedness holds automatically as long as dominance is satisfied. \end{rem}       

\begin{rem}[\textbf{Domination and variable selection consistency}] Even if there exist algorithms for $L$ resp. $\tilde L$ that both have the variable selection consistency property (see \cite{bu}), there would be no guarantee that there would be no empirical singular parts in real applications, i.e., \textbf{for finite} $n$. \end{rem}  

\begin{rem}[\textbf{Surrogate losses and singular parts}] Although the usage of surrogate loss functions is common for example in classification settings, we are not sure if a ''wrong'' surrogate may even lead to singular parts of the column measure corresponding to the original loss (e.g. the $0/1-$loss) w.r.t. the column measure corresponding to the surrogate loss.  \end{rem}

\begin{ex}\label{stabselex} For the moment, ignore singular parts and consider a simple example where we apply the Stability Selection of \cite{hofner15} implemented in $\mathsf{R}-$package \texttt{stabs} (\cite{stabs}) to Boosting models with different loss functions. \ \\

More specifically, consider the following four losses $L_2$ and $L_1$ for squared and absolute loss respectively, $L_{\tau}$ for check loss and $L_{huber}$ for Huber loss. Note that all four losses agree that an error-free prediction would be optimal. Then, in a simulated data set\footnote{Details and $\mathsf{R}-$code are available from the author upon request.} with $n=100$ observations from a Gaussian linear model with $p=50$ metric (Gaussian) predictors, $s^0=10$ true non-zero coefficients and signal-to-noise ratio 2, we obtain no consensus concerning a stable model for these four loss functions. As a consequence, we conclude that although all four losses head for ''optimal'' prediction, there is no canonical candidate loss $L$ which could be chosen as reference for our hard ranking loss $\tilde L$. In particular, there is no guarantee that $L_2-$Boosting selects all $\tilde L-$relevant variables.  \end{ex} \

Summing up. when trying to solve the problem of sparse empirical minimization of loss tilde L by a (generalized) rejection scheme starting from loss L, we are led to the following central questions: \ \\

\textit{\textbf{(Q1) How can we control for singular parts of the (unknown) column measure $\nu^{(\tilde L)}$ w.r.t. the (unknown) column measure $\nu^{(L)}$?}}
\textit{\textbf{(Q2) How can we get hand on $\nu^{(\tilde L)}$ at all in order to warrant that we change measure $\nu^{(L)}$ to the correct one $\nu^{(\tilde L)}$?}}

\section{SingBoost: Boosting with singular parts for any target loss} \label{singboostsec} \

The relevance of a certain predictor in a Boosting model w.r.t. a loss $L$ is simply determined by dividing the number of iterations in which this variable has been selected by the number of Boosting iterations. As we already discussed, each Boosting iteration is a rejection sampling step itself, comparing weak learners w.r.t. each predictor and selecting the one that performs best, evaluated in the loss function $L$. Since trying to control for singular parts on the level of whole Boosting models is inappropriate in the presence of singular parts, we need to interfere deeper into the Boosting algorithm itself and force that $\tilde L$ is already respected there. \ \\

Settling questions (Q1) and (Q2) amounts to two steps, i.e., (a) a strategy warranting that we visit singular parts of $\nu^{(\tilde L)}$ w.r.t. $\nu^{(L)}$ sufficiently often and (b) the actual change of measure from $\nu^{(L)}$ to $\nu^{(\tilde L)}$. \ \\

\subsection{Change of measure by secant steps} \

To this end it helps to recall how $L_2-$Boosting can be seen as a functional gradient decent procedure. We again refer to \cite{zhao04} for the following lemma.    

\begin{lem} The column selection done in each $L_2-$Boosting-step by finding the predictor column maximizing the absolute value of correlation to the current response variable is exactly selecting the column which for given step length $\kappa$ achieves the largest descent in the loss function $L$.     \end{lem} \

For column selection according to more general losses $\tilde L$, one now has two options: Either one can compute the gradient for each candidate column and finds the one with the largest gradient, or if computation of the gradient is either too burdensome, or, moreover if the gradient does not even exist, we can evaluate the secant, i.e., $\Delta \tilde L$ and select the column with steepest descent. This provides a strategy for (Q2).             \ \\

In the sequel we confine ourselves to this latter strategy as this in particular covers the loss function $L_n^{hard}$ (Equation (\ref{rankparopt})) from the continuous ranking problem.  \ \\

Still, we note that in principle all strategies to speed up the computation of $\Delta \tilde L$ are eligible here. 

As evaluation of $\Delta \tilde L$ tends to be more expensive than computation of correlations (as needed for $L_2-$Boosting), for some suitably chosen value $M$ (e.g., 5, 10), we only perform this every $M-$th iteration.    \ \\

\subsection{Coping with singular parts} \

As a convenient side effect, when checking all columns in the search for the best candidate column, we automatically take care about the singular parts. \ \\

Heuristically, while keeping the efficient structure of componentwise $L_2-$Boosting, we are getting a chance to have variables being detected that correspond to a singular part of $\nu^{(\tilde L)}$ w.r.t $\nu^{(L)}$ being detected.  

Even more, due to the fact that classical Gradient Boosting is greedy and never drops variables once they are selected, one may even organize this in Algorithm \ref{singboost}, providing a constructive answer to (Q1). \ \\

\begin{algorithm}[H] 
\label{singboost}
\textbf{Initialization:} Data $\D^{sing}$, step size $\kappa \in ]0,1]$, number $m_{iter}$ of iterations, number $M \le m_{iter}$ (each $M-$th iteration is a singular iteration), target loss $\tilde L$ (as part of a \texttt{family} object \texttt{singfamily} in the infrastructure of $\mathsf{R}-$package \texttt{mboost}, \cite{mboost}) and binary variable \texttt{LS}\;
Set $$ \runs=\left\lfloor \frac{m_{iter}}{M} \right\rfloor \; $$ 
Set $f^{(0)}:=0$\;
\For{$k=1,...,\runs$}{
\eIf{\texttt{LS==FALSE}}{Perform a single step of Gradient Boosting w.r.t. $\tilde L$ on the residuals of model $\hat f^{((k-1)M)}$\; }{Evaluate all simple least squares fits on the residuals of model $\hat f^{((k-1)M)}$ w.r.t. $\tilde L$ and take the best one\;}
Get the weak model $\hat g^{((k-1)M+1)}$ and update the model via \begin{center} $ \displaystyle \hat f^{((k-1)M+1)}=\hat f^{((k-1)M)}+\kappa \hat g^{((k-1)M+1)} $\; \end{center} 
Perform $(M-1)$ steps of $L_2-$Boosting starting with the residuals w.r.t. the model $\hat f^{((k-1)M+1)}$\; 
Get the updated model $\hat f^{(kM)}$
}
\caption{SingBoost}\index{SingBoost!Pseudocode}
\end{algorithm} \ \\ 

\begin{rem} The argument \texttt{LS=FALSE} which performs a single Gradient Boosting step w.r.t. $\tilde L$ is only meaningful in simulation settings if one wants to compare the performance of different baselearners. In general (and when we actually need SingBoost), $\tilde L$ does not have a corresponding ''pure'' Boosting algorithm. Then setting \texttt{LS=TRUE} uses simple least squares baselearners in the singular steps, but compares them according to $\tilde L$. \end{rem}

\begin{rem} While as already noted loss $L^{hard}$ from Equation (\ref{rankparopt}) is not even continuous (and nor are typically losses for localized ranking), they can at least be computed with  $\mathcal{O}(n\ln(n))$ steps thanks to a lemma from the first author's review paper (\cite{TW19b}). \end{rem} 

\begin{rem}[\textbf{Why $L_2-$Boosting?}] \label{whyl2rem} The question whether SingBoost leads to an effective variable selection when using the column measure $\hat \nu^{(L)}$ from $L_2-$Boosting can be answered in the same manner as the question of effectivity of rejection sampling, i.e., the closer $\nu^{(L)}$ is w.r.t. $\nu^{(\tilde L)}$ the more effective the procedure; like in rejection sampling this could even be quantified in terms of mutual Kullback-Leibler information; in case of the continuous ranking problem, there is evidence that indeed $L_2-$Boosting should provide a good candidate $\nu^{(L)}$, as a perfect $L_2-$regression fit would also be perfect for the continuous ranking problem (see \cite{clem18}).  \end{rem} 

\begin{rem}[\textbf{Sparse Boosting}]\index{Boosting!Sparse Boosting} Although not intended to detect potential singular parts by tailoring the model selection to another loss function, the Sparse Boosting algorithm (\cite{bu06a}) can be regarded as related work. The aim of this algorithm is to provide sparser models by trying to minimize the out-of-sample prediction error. Since this is not directly possible, they use the degrees of freedom defined as $\tr(\mathcal{B}_m)$ for the Boosting operator \begin{center} $ \displaystyle \mathcal{B}_m=I-(I-\kappa \mathcal{H}_{\hat j_m})\hdots(I-\kappa \mathcal{H}_{\hat j_1})$ \end{center} to estimate the model complexity where $\mathcal{H}_{\hat j_k}$ is the hat matrix corresponding to the weak model which is based on variable $\hat j_k$ (see also \cite{bu03} for further details on the Boosting operator). In fact, they choose the variable that minimizes the residual sum of squares, penalized by the model complexity, and then proceed as in $L_2-$Boosting. As they pointed out, this leads to \textbf{possibly different choices of the best variable} in each iteration which in some sense can also be understood as some kind of SingBoost with a \textbf{sequence of target loss functions} $\tilde L^{(k)}$ that equal their penalized criterion in every iteration $k=1,...,m$. \end{rem}  \

Let us highlight again the most important algorithmic fact concerning the loss function $\tilde L$. The only property of $\tilde L$ that we require to get a computationally tractable algorithm is that $\tilde L$ is easy to evaluate. Apart from this, no (local) differentiability property, not even continuity of $\tilde L$ is required. In other words, SingBoost extends the Gradient Boosting setting that was restricted to convex and almost-everywhere differentiable loss functions to the situation of (nearly) arbitrary loss functions, i.e., we can think of SingBoost as a \textbf{gradient-free Gradient Boosting algorithm}. \ 

\begin{rem} In order not to get trapped by pathological jumps in $\tilde L$, the loss should be ''asymptotically smooth'' in the sense that the aggregation induced by the summation in the passage to the risk should warrant that the risk be smooth in the predictions $\hat Y$. This is true for the hard continuous ranking loss as soon as the underlying distributions are absolutely continuous. \end{rem} \

As for the updating step, for any loss function $\tilde L$ which has a structure such that $\tilde L(|r|)$ is monotonically increasing, we can transfer the selection criterion from componentwise least squares Boosting to SingBoost by taking \begin{center} $ \displaystyle \hat j_k=\argmin_j\left(\sum_i \tilde L(r_i^{(k-1)}-X_{i,j}\hat \beta(j))\right) .$ \end{center} 

Evidently, we may even cover more general loss functions $\tilde L$ of form $\tilde L(Y,\hat Y)$ and a more general learning rate $\kappa$ by selecting $\hat j_k$ by \begin{equation} \label{varselrank} \hat j_k=\argmin_j\left(\sum_i \tilde L(Y_i,X_i\hat \beta^{(k-1)}+\kappa \hat \beta(j))\right), \end{equation} and consequently obtain a strong model like in general Boosting algorithms. In fact, this criterion can be seen as a secant approximation to $\nabla \tilde L$. \ \\

As for the complexity of SingBoost, see the following lemma. \ \\

\begin{lem}[\textbf{Complexity of SingBoost}]\index{SingBoost!Complexity}\label{singcomplexlem} Assume that the evaluation of the target loss function $\tilde L$ requires $\mathcal{O}(c_n)$ operation for some $c_n \ge n \ \forall n$. Then the complexity of SingBoost is $\mathcal{O}(mc_np)$. \end{lem}

\begin{rem}\label{singcomplex} Since $L_2-$Boosting is of complexity $\mathcal{O}(mnp)$, the loss in performance (assuming an equally excellent implementation) for SingBoost compared to $L_2-$Boosting gets negligible if $c_n=\mathcal{O}(n)$. \\
In the case of the hard ranking loss, we have $c_n=\mathcal{O}(n\ln(n))$ due to the lemma from the first author's review paper \cite{TW19b}, so for real-data applications where $n$ is usually rather small, the complexity of SingBoost is comparable with that of $L_2-$Boosting, in particular as not each iteration is a singular iteration. \end{rem}   

\begin{rem}[\textbf{Step size}] At this stage it is not clear \`{a} priori how to select step length $\kappa$. In fact, it is not obvious that guidelines from classical Gradient Boosting (see \cite{bu07}) can be taken over without change. Preliminary evaluations show that at least for the continuous ranking problem the guidelines for $L_2-$Boosting seem to be reasonable for SingBoost, too, and we delegate a more thorough discussion of this choice to further research.   \end{rem}

\subsection{Connection to resampling techniques and rejection sampling} \

As announced, having explicated the two measures $\nu^{(\tilde L)}$ and $\nu^{(L)}$ and the task of a change of measure respecting singular parts, this opens the door to many well-established simulational techniques to achieve this task, such as resampling and rejection sampling. \ \\

As an example, we spell out a variant of the classical rejection sampling  (e.g. \cite[Alg. 3.4]{ripley}) which could be used for this task. \ \\

Assume that through some initialization phase where one uses SingBoost as defined in Algorithm \ref{singboost} we have estimates $\hat \nu^{(\tilde L)}$ and $\hat \nu^{(L)}$ for (the counting densities of)  $\nu^{(\tilde L)}$ and $\nu^{(L)}$, and derived from these the index sets  \begin{center} $ \displaystyle J_c:=\{j \in \N \ | \ \hat \nu^L(\{j\})>0\}, \ \ \ J_s:=\{1,...,p\}\setminus J_c. $ \end{center} Based on these index sets, we define the masses \begin{center} $ \displaystyle w_c:=\frac{\# J_c}{p}, \ \ \ w_s:=1-w_c , $ \end{center} \begin{center} $ \displaystyle  W_s:=\sum_{j\in J_s} \nu^{(\tilde L)}(\{j\}) $ \end{center} and the bound \begin{center} $ \displaystyle H:=\max_{j \in J}(\hat \nu^{(\tilde L)}(\{j\}) /\hat \nu^{(L)}(\{j\})) $. \end{center} \ 

Then a singular-part-adjusted variant of rejection sampling could be implemented as in Algorithm \ref{RejSampling}.

\begin{algorithm}
\label{RejSampling}
\textbf{Initialization:} Active candidate set $J^\natural=\{1,...,p\}$\;
\textbf{(1)} Sampling of a variable $B \sim \text{Binom}(1,w_c)$; according to $B$ switch between cases (a)  (for $B=0$) and (b) (else) \;
\textbf{(2a)} Sampling an index $j_0$ from $J_s$ with column weights $\hat \nu^{(\tilde L)}(\{j\})/W_s$\;
\textbf{(3a)} Usage of $j_0$ as new column\;
\textbf{(4a)} Set $J^\natural=\{1,...,p\}$\; 
       Jump to (1)\;
\textbf{(2b)} Selection of the next candidate column $j_0$ by $\hat \nu^{(L)}|_{J^{\natural}}$ resp. $L_2-$Boosting (Algorithm \ref{l2Boosting}) among all columns in $J^\natural$\;
\textbf{(3b)} Computation of $h(j_0):=\hat \nu^{(L)}(j_0) / \hat \nu^{(\tilde L)}(j_0)/H$\;
\textbf{(4b)} Sample $U \sim Unif([0,1])$\;
\textbf{(5b)} \eIf{$U \le h(j_0)$}{Accept $j_0$ from (2b) as column and set $J^\natural=\{1,...,p\}$}{Set $J^\natural=J^\natural\setminus \{j_0\}$\;}
\textbf{(6b)} Jump to (1) \\

{\footnotesize In order to give all columns in $J_s$ a chance to be visited, $\hat \nu^{(\tilde L)}$ must put a minimal positive mass on each $j\in J_s$.}

\caption{Rejection Sampling from column measures}
\end{algorithm}

\section{Asymptotic properties of SingBoost} \label{asysing} \

Our modifications of the standard $L_2-$Boosting do not affect the theoretical properties derived for the original $L_2-$Boosting. More precisely, we can translate the theorems \cite[Thm. 1]{bu06} and \cite[Thm. 12.2]{bu} on estimation and prediction consistency by adding a so-called \textbf{Corr-min condition}. \ \\

\cite[Thm. 1]{bu06} use the following scheme going back to Temlyakov (\cite{temlyakov}): For \begin{center} $ \displaystyle \langle f,g \rangle_{(n)}:=\frac{1}{n}\sum_i f(X_i)g(X_i), $ \end{center} i.e., an empirical version of the inner product $\langle f,g \rangle$ on the space $L_2$. B\"{u}hlmann identifies $L_2-$Boosting as the iterative scheme \begin{equation} \label{temlyakov} \begin{gathered} \hat R_n^0f=f, \ \ \ \hat R_n^1f=f-\langle Y,g_{\widehat{j_1}} \rangle_{(n)}g_{\widehat{j_1}}, \\ \hat R_n^mf=\hat R_n^{m-1}f-\langle \hat R_n^{m-1}f,g_{\widehat{j_m}} \rangle_{(n)}g_{\widehat{j_m}} \ \forall m \ge 2, \end{gathered} \end{equation} where $g_j$ represents the baselearner w.r.t. variable $j$ (it is assumed that $\langle g_j,g_j \rangle_{(n)}=1$ for all $j$) and $\hat R_n^m$ is the empirical residual after the $m-$th iteration.  \ \\

For the special case of $L_2-$Boosting, the selected variables are, as described in \cite{bu06}, \begin{equation}\label{temlyakovvarsel} \widehat{j_1}=\argmax_j(|\langle Y,g_j \rangle_{(n)}|), \ \ \ \widehat{j_m}=\argmax_j(|\langle \hat R_n^{m-1}f,g_j \rangle_{(n)}|) \ \forall m \ge 2. \end{equation} 

So, the only difference between $L_2-$Boosting in SingBoost in the language of the Temlyakov scheme is the variable selection in the singular iterations, leading to other residuals according to (\ref{temlyakov}). We mimic \cite[Thm. 1]{bu06} and propose the following result for \textbf{random design of the regressor matrix}: \

\begin{thm}[\textbf{Estimation consistency of SingBoost}] \label{asyconssing} Let us define the model \begin{center} $ \displaystyle Y_i=f_n(X_i)+\epsilon_i=\sum_{j=1}^{p_n} \beta_{n,j}X_{i,j}+\epsilon_i $ \end{center} for $X_i \in \R^{p_n}$ i.id. with $\E[||X_{\cdot,j}||^2]=1$ for all $j=1,...,p_n$ and error terms $\epsilon_1,...,\epsilon_n$ that are i.id., independent from all $X_i$ and mean-centered. Let $X$ denote a new observation which follows the same distribution as the $X_i$, independently from all $X_i$. Assume (E1)-(E4) from \cite[Thm. 1]{bu06}, i.e.,  \ \\

\textbf{(E1)} $\exists \xi \in ]0,1[, C \in ]0,\infty[$ such that $p_n=\mathcal{O}(\exp(Cn^{1-\xi}))$ for $n \rightarrow \infty$, \ \\

\textbf{(E2)} $\sup_n\left(\sum_{j=1}^{p_n} |\beta_{n,j}|\right)<\infty$, \ \\

\textbf{(E3)} $\sup_{j,n}(||X_{\cdot,j}||_{\infty})<\infty$, \ \\

\textbf{(E4)} $\E[||\epsilon||^{\delta}]<\infty$ for $\delta>\frac{4}{\xi}$ with $\xi$ from (E1). \ \\

And in addition assume \ \\

\textbf{(E5)} $\exists a>0 \ \forall m,n \ \exists \tilde C \ge a:$ \begin{equation} \label{corrmin} |\langle \hat R_n^{m-1}f,g_{\widehat{j_m}}\rangle_{(n)}| \ge \tilde C\sup_j(|\langle \hat R_n^{m-1}f,g_j \rangle_{(n)}|). \end{equation} 

Then, denoting the SingBoost model after the $m-$th iteration based on $n$ observations by $\hat f_n^{(m)}$, it holds that \begin{center} $ \displaystyle \E_X[||\hat f_n^{m_n}(X)-f_n(X)||^2]=o_P(n^0) $ \end{center} for $n \rightarrow \infty$ provided that $(m_n)_n$ satisfies that $m_n \rightarrow \infty$ sufficiently slowly for $n \rightarrow \infty$. \end{thm}

\begin{rem}[\textbf{Step size}]\index{Boosting!$L_2-$!Step size}\index{SingBoost!Step size} The proof implicitly assumed the learning rate $\kappa=1$. Of course, we can follow the same lines as in \cite[Sec. 6.3]{bu06} to achieve the main result for learning rates $\kappa<1$. \end{rem}  

\begin{rem} We did not explicitly use the fact that only each $M-$th iteration is a singular iteration, i.e., for all other $m$, nothing changes compared to \cite[Thm. 1]{bu06}. Indeed, asymptotically, the statement of the theorem does not change, only the bound for finite $n$ would be affected, but we do not see an advantage in a tedious case-by-case-analysis. \end{rem} \

The constant $a>0$ from the Corr-min condition (E5) was absorbed during the proof when one forced the sequence $(m_n)_n$ to grow sufficiently slowly to get an upper bound of order $o_P(n^0)$ for $||\tilde R_n^mf||$. However, in the following theorem that is based on \cite[Thm 12.2]{bu}, a bit more work is necessary to adapt it to SingBoost, directly appearing in the convergence rate at the end. \ 

\begin{thm}[\textbf{Prediction consistency of SingBoost}]\index{SingBoost!Prediction consistency} \label{asypredsing} Let us define the model \begin{center} $ \displaystyle Y_i=f_n(X_{n;i})+\epsilon_i=\sum_{j=1}^{p_n} \beta_{n,j}X_{n;i,j}+\epsilon_i $ \end{center} for \textbf{fixed design} of the regressor matrix and error terms $\epsilon_1,...,\epsilon_n$ that are i.id., independent from all $X_i$ and mean-centered. Let $X$ denote a new observation which is independent from all $X_i$. Assume conditions (P1)-(P4) from \cite[Thm. 12.2]{bu}, i.e., \ \\

\textbf{(P1)} The number $p_n$ satisfies $\frac{\ln(p_n)}{n} \rightarrow 0$ for $n \rightarrow \infty$, \ \\

\textbf{(P2)} The true coefficient vector is sparse in terms of the $l_1-$norm, i.e., \begin{center} $ \displaystyle ||\beta_n||_1=\sum_j |\beta_{n,j}|=o\left(\sqrt{\frac{n}{\ln(p_n)}}\right) $ \end{center} for $n \rightarrow \infty$, \ \\

\textbf{(P3)} It holds that \begin{center} $ \displaystyle \frac{1}{n} \sum_i X_{n;i,j}^2=1 \ \forall j=1,...,p_n $ \end{center} and \begin{center} $ \displaystyle \frac{||X\beta||_2^2}{n}=\frac{1}{n}\sum_i f_n^2(X_{n;i}) \le C<\infty $ \end{center} for all $n \in \N$, \ \\

\textbf{(P4)} The errors are i.id. $\mathcal{N}(0,\sigma^2)-$distributed for all $n \in \N$ \ \\

and in addition, assume the Corr-min condition (E5) of the previous theorem. Then for \begin{center} $ \displaystyle m_n \rightarrow \infty, \ \ \ m_n=o\left(\sqrt{\frac{n}{\ln(p_n)}}\right) $ \end{center} for $n \rightarrow \infty$, it holds that \begin{center} $ \displaystyle \frac{||X(\hat \beta_n^{(m_n)}-\beta_n)||_2^2}{n}=o_P(n^0) $ \end{center} as $n \rightarrow \infty$ where $\hat \beta_n^{(m_n)}$ is the coefficient vector based on $n$ observations and after the $m_n-$th SingBoost iteration. \end{thm}   \ \\

Both proofs are delegated to Section \ref{techproofs}. \

\begin{rem}[\textbf{Discussion of the Corr-min condition (E5)}] Imposing the Corr-min condition could be seen as restrictive; in applications we studied so far, we have not yet encountered instances where we could not satisfy this condition though. In theory one could think of good predictive models found by a Gradient Boosting algorithm for $\tilde L$ where the Corr-min conditions is violated.  \ \\

This is already true if one concerns quantile Boosting. If one uses the implementation of quantile Boosting in the $\mathsf{R}-$package \texttt{mboost}\index{Boosting!\texttt{mboost}} by using either \texttt{QuantReg()} or, in the case of $L_1-$Boosting, \texttt{Laplace()} as family object, then it may happen that only the intercept is being selected in some (or even in all) iterations. Therefore, if we based SingBoost on quantile Boosting, the Corr-min condition will always be violated as soon as one only selects the intercept column. \ \\

This situation is excluded in our Singboost algorithm, though: It always uses simple least squares models as baselearners which are directly based on the correlation of the current residual with each single variable. Therefore, if we have variables that are (nearly) perfectly uncorrelated with the current residual, then the resulting coefficient would be (in a very close neighborhood of) zero. Although we cannot exclude this case, it would be very unlikely that such a quasi-zero-baselearner would pass the rejection step and enter the SingBoost model. \index{SingBoost!Corr-min condition|)} \end{rem}  

\begin{rem}[\textbf{Variable selection consistency}] In contrast to the Lasso (\cite{tibsh96}) for which there exist results on variable selection consistency requiring a beta-min condition and an irrepresentability condition (see \cite{bu}), \cite{vogt} recently showed that we may not expect a corresponding result for $L_2-$Boosting, even if the restricted nullspace property is satisfied, i.e., the intersection set of the nullspace of $X$ and the cone \vspace{0.2cm} \begin{center} $ \displaystyle C(S,L):=\{\beta \in \R^p \ | \ ||\beta_{(S^0)^c}||_1 \le L||\beta_{S^0}||_1\} $ \end{center} \vspace{0.2cm}  is just the element $0_p$, it is, \textbf{$L_2-$Boosting is not variable selection consistent}\index{Boosting!$L_2-$!Variable selection inconsistency}. The authors provide an example for a regressor matrix with the restricted nullspace condition but where $L_2-$Boosting reliably selects the wrong columns. \ \\

As in SingBoost through the singular iterations we do include columns which are never selected in $L_2-$Boosting we cannot tell whether the example in \cite{vogt} has a translation to SingBoost.  \end{rem} \ \\

\section{Coefficient paths} \

Implementations of sparse learning algorithms like the Lasso or the elastic net are capable to supply the user with a regularization path of the coefficients, like in the $\mathsf{R}-$package \texttt{glmnet} (\cite{friedman10}). Those paths show the evolution of the coefficients for different values of the regularization parameter $\lambda$, more precisely, in the cited package, one can (amongst other options) let the values of the non-zero coefficients be plotted against the natural logarithm of the sequence of regularization parameters (in an ascending order, i.e., the paths evolve from the right to the left).  \ \\

For Boosting, similar coefficient paths are available in the $\mathsf{R}-$package \texttt{mboost}\index{Boosting!\texttt{mboost}}\index{Boosting!Coefficient paths|(}. Since Boosting does not incorporate a regularization parameter, those Boosting coefficient paths are given by the values of the (non-zero) coefficients vs. the current Boosting iteration. Clearly, since only one variable, not counting the intercept, changes in one Boosting iteration, the coefficient paths, except one, are flat when jumping from on iteration to the next. \ \\

These coefficient paths are a valuable diagnostic tool providing a visualization of selection order, and inspecting the flattening yields a indication for stopping. \ \\

\begin{figure}
\centering
\includegraphics[width=7cm]{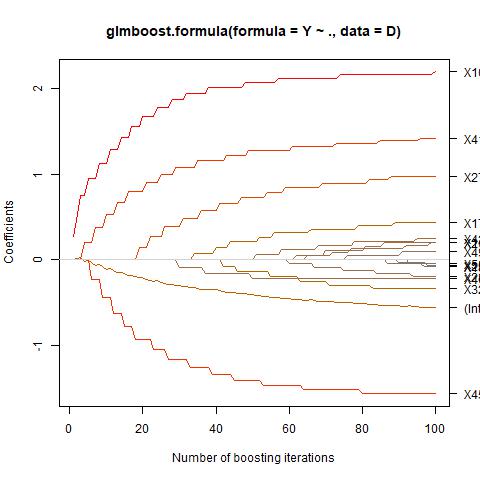}
\caption[Coefficient paths provided by \texttt{glmboost}]{Coefficient paths for $L_2-$Boosting} \index{Boosting!Coefficient paths|)}     \label{glmpaths}
\end{figure}

\begin{figure}
\centering
\includegraphics[width=7cm]{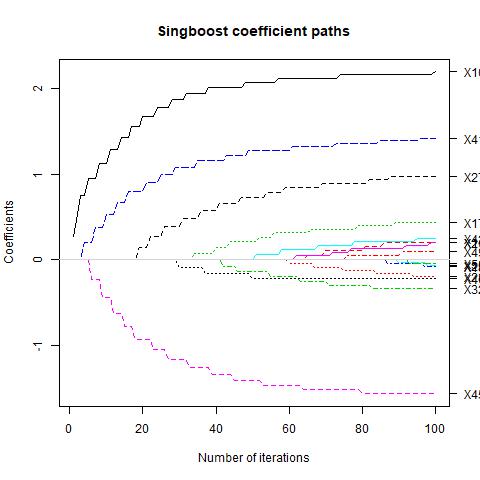}
\caption[Coefficient paths provided by \texttt{glmboost}]{Coefficient paths for SingBoost}   \label{singpaths}
\end{figure} 

As our algorithm SingBoost only replaces the selection mechanism of a single column but otherwise maintains the structure of a general Functional Gradient Boosting algorithm, these coefficient paths are also available for SingBoost. \ \\

The coefficient paths in Figure \ref{glmpaths} have been drawn in $\mathsf{R}$, \cite{Rlang}, using the \texttt{plot} command in the $\mathsf{R}-$package \texttt{mboost}. The coefficient paths in Figure \ref{singpaths} were generated by our function \texttt{singboost.plot} on the same data set and with the same configurations as in \texttt{glmboost} (especially setting \texttt{singfamily=Gaussian()}). In the simulation, we used a data set with $n=100$ observations from a Gaussian linear model with $p=50$ metric (Gaussian) predictors, $s^0=10$ true non-zero coefficients and signal-to-noise ratio 2. \ \\ \ \\

\section{Outlook} \

So far, we developed the algorithm SingBoost enhancing the original $L_2-$Boosting algorithm such that it may select variables that $L_2-$Boosting does ignore. On the other side, we already mentioned that apart from the actual column selection SingBoost respects the structure of a generic Functional Gradient algorithm. This said, it is not surprising that SingBoost also inherits the proneness to overfitting and the instability from generic Gradient Boosting. \ \\

For generic Gradient Boosting, these issues have been addressed through Stability Selection (\cite{bu10}, \cite{hofner15}). It is obvious that a similar, but adapted strategy should also be helpful in context of SingBoost. \ \\

In addition, if there are singular parts, the corresponding coefficients can only be updated in the singular steps and are therefore disadvantaged, depending on the frequency $M$ of singular steps. \ \\

Indeed, the present authors have set up an adapted Stability concept for this purpose. This concept though is valid in a more general learning framework, not only restricted to Gradient Boosting, hence a thorough discussion of this concept would be out of scope for this paper. Instead, we rather refer to \cite{TWphd} and a subsequent paper. \ \\

In particular, we decided not to provide simulation evidence for our SingBoost algorithm alone and instead rather to focus on the theoretical properties of SingBoost and the more general column measure framework here, as without Stability criterion the empirical results obtainable from mere SingBoost could be misleading. \ \\ 

To include the main idea of our adaptation of the Stability Selection, note that the original Stability Selection as introduced by \cite{bu10} amounts to a clever application of re-/subsampling ideas where only those columns survive which are selected in the majority of resampling instances. Our stabilization concept also uses subsampling but in addition may also attach weights to the resampled instances according to their $\tilde L-$performance. \ \\ \\

\section{Conclusion} \

Motivated by the issue that different loss functions lead to different variable selection when performing empirical risk minimization, we introduced measures on the rows resp. the columns of a data matrix and called them row measure resp. column measure to reflect the in-/exclusion of certain columns (and rows) into the model in a mathematically coherent way. It turns out that many traditional learning procedures may be cast into this unified framework; as examples we combine the One-Step estimators with column measures giving the Expected One-Step. \ \\

The column measure framework also renders explicit issues in model selection when performance is measured according to loss $\tilde L$ but selection is done according to loss $L$. These issues can be identified mathematically as singular parts of column measure $\nu^{(\tilde L)}$ w.r.t. to $\nu^{(L)}$.   \ \\

Identifying each iteration of componentwise Boosting as rejection step itself, we proposed the algorithm SingBoost that includes singular steps where a linear baselearner is evaluated in the target loss $\tilde L$ so that we get the chance to select variables from potential singular parts. These additional singular steps can cope with non-differentiable, even non-continuous losses $\tilde L$, hence cover in particular the hard continuous ranking loss from Equation (\ref{rankparopt}). \ \\

Encouragingly, we could show that our new SingBoost algorithm enjoys the same attractive statistical asymptotic properties as to prediction and estimation consistency as $L_2-$Boosting.  \ \\

\subsection{Acknowledgements} \

Most of this paper is part of T. Werner's PhD thesis at Oldenburg University under the supervision of P. Ruckdeschel. \ \\  \ \\

\section{Technical proofs} \label{techproofs} \

\begin{bew}[Theorem \ref{asyconssing}] The proof follows the same steps as the proof of \cite[Thm. 1]{bu06}. For his proof, B\"{u}hlmann uses two Lemmata (\cite[Lemma 1+2]{bu06}). His Lemma~1 does not include variable selection and indeed holds for our case. His Lemma~2 also holds since an upper bound of the expression $\langle g_{\widehat{j_m}}, g_j \rangle_{(n)}$ is used which is independent of $\widehat{j_m}$. \ \\

Using (E5), Equation (6.13) in \cite{bu06} changes to \begin{equation} \label{613new} |\langle \tilde R_n^mf,g_{\widehat{j_{m+1}}}\rangle | \ge \tilde C \sup_j(|\langle \tilde R_n^mf,g_j \rangle |)-(1+\tilde C)(2.5)^m\zeta_n C_* \end{equation} with $\zeta_n$, $C_*$ from \cite[Lem. 2]{bu06}, on the set $A_n:=\{\omega \ | \ |\zeta_n(\omega)|<0.5 \}$. Now, we have to consider an analog of the set $B_n$ defined on p. 580 in \cite{bu06}, namely \begin{equation} \label{Bnnew} \tilde B_n:=\{\omega \ | \ \sup_j(|\langle \tilde R_n^mf,g_j \rangle |)>2(\tilde C)^{-1}(1+\tilde C)(2.5)^m\zeta_n C_*\}, \end{equation} and can conclude that, using (\ref{613new}), it holds that \begin{center} $ \displaystyle |\langle \tilde R_n^mf,g_{\widehat{j_m}} \rangle | \ge 0.5\tilde C \sup_j(|\langle \tilde R_n^mf,g_j \rangle|) $ \end{center} on $A_n \cap \tilde B_n$. That means, we have $b=0.5\tilde C \ge 0.5a>0$ in Equation (6.2) of \cite{bu06}. This is less than the constant $0.5$ established in \cite{bu06}, but it does not matter as long as it is bounded away from zero. Then, we get an analog of Inequality (6.15) of \cite{bu06} which is \begin{center} $ \displaystyle ||\tilde R_n^mf||\le B(1+mb^2)^{-b/2(2+b)}=o(n^0) $ \end{center} for $n \rightarrow \infty$ provided that $m_n \rightarrow \infty$ as assumed.  \ \\

On $\tilde B_n^c$, we proceed in the same manner as \cite{bu06} since $\tilde B_n$ can be identified with the set notated as $\mathcal{A}(\mathcal{D},b_m)$ in \cite{temlyakov} which supplies a recursive scheme, leading to the same upper bound as in Inequality (6.16) in \cite{bu06} since the variable selection is absorbed when applying the ''norm-reducing property'' given in Equality (6.3) in \cite{bu06}. \ \\

The rest of the proof follows exactly the same steps as the proof in \cite{bu06} since $\widehat{j_m}$ only appears again when bounding the quantity $A_n(m)=||\hat R_n^mf-\tilde R_n^mf||$, but $||g_{\widehat{j_m}}||$ is bounded by one and \cite[Lem. 2]{bu06} holds anyway, yielding the last inequality on page 580 in \cite{bu06}. \begin{flushright} $_\Box $ \end{flushright} \end{bew} \ 

Note that the equality in the second last display on page 579 in \cite{bu06} contains a typo in the index. More specifically, in Inequality (6.13) in \cite{bu06}, one should replace $\widehat{j_m}$ by $\widehat{j_{m+1}}$ and correspondingly in Inequality (6.14). \ \\ 

\begin{bew}[Theorem \ref{asypredsing}] We follow the same steps as in the proof of \cite[Thm. 12.2]{bu} which again is based on a Temlyakov scheme as before. First, \cite[Lem. 12.1]{bu} needs to be modified to: \ \\

\begin{lem} If there exists $0<\phi<0.5$ such that \begin{center} $ \displaystyle \max_j(|\langle \hat R^{m-1}f,X_{\cdot,j}\rangle_{(n)} |) \ge 2\Delta_n \phi^{-1} (\tilde C)^{-1}, $ \end{center} then it holds that \begin{center} $ \displaystyle |\langle \hat R^{m-1}f,X_{\cdot,\widehat{j_m}}\rangle_{(n)} | \ge \tilde C(1-\phi)\max_j(|\langle \hat R^{m-1}f,X_{\cdot,j} \rangle_{(n)}|), $ \end{center} \begin{center} $ \displaystyle |\langle Y-\hat f^{(m-1)},X_{\cdot,\widehat{j_m}} \rangle_{(n)}| \ge \tilde C(1-\phi/2)\max_j(|\langle \hat R^{m-1}f,X_{\cdot,j} \rangle_{(n)} |). $ \end{center} \end{lem} 

The proof just modifies the proof of \cite[Lem. 12.1]{bu} at obvious places. \ \\

The quantities $a_k$ and $d_k$ defined on page 422 in \cite{bu} indeed depend on the selected column, but Inequality (12.28) in \cite{bu} holds as well and since the relationship of $a_k$ and $d_k$ does not change with the concrete column selection and since both are bounded from above at the end, the singular iterations do not affect the proof structure when facing $a_k$ and $d_k$. \ \\

\cite[Lem. 12.2]{bu} has to be modified to: 

\begin{lem} \label{122surro} If there exists $0<\phi<1/2$ such that \begin{center} $ \displaystyle \max_j(|\langle \hat R^{m-1}f,X_{\cdot,j}\rangle_{(n)} |) \ge 2\Delta_n \phi^{-1}(1-\phi/2)^{-1}(\tilde C)^{-1}, $ \end{center} then it holds that \begin{center} $ \displaystyle a_m \le a_{m-1}-(1-\phi)d_m^2. $ \end{center} \end{lem} 

The proof also follows the same lines as the proof of \cite[Lem. 12.2]{bu}. \ \\

For shortness, we only specify the necessary modifications here. \ \\

Inequality (12.33) in \cite{bu} changes to \begin{center} $ \displaystyle d_k \ge \frac{(1-\phi/2)a_{k-1}\tilde C}{b_{k-1}} $ \end{center} and instead of Inequality (12.34) in \cite{bu}, we get \begin{center} $ \displaystyle a_kb_k^{-2} \le a_{k-1}b_{k-1}^{-2}(1-C_{\phi}^2a_{k-1}b_{k-1}^{-2}\tilde C) $ \end{center} with $C_{\phi}$ as in \cite{bu}.  Instead of $B_n(m)$ as given in Equality (12.36), we define \begin{center} $ \displaystyle \tilde B_n(m):=\bigcap_{k=1}^m \{\max_j(|\langle \hat R^{k-1}f^0,X_{\cdot,j} \rangle_{(n)} |) \ge 2(\tilde C)^{-1}\phi^{-1}(1-\phi/2)^{-1} \Delta_n \} $ \end{center} and, analogously to Inequality (12.37) in \cite{bu}, we conclude that \begin{center} $ \displaystyle a_mb_m^{-2} \le ||\beta_n^0||_1^{-2} (1+C_{\phi}^2m\tilde C)^{-1}. $ \end{center}  

Inequality (12.38) in \cite{bu} is modified to \begin{center} $ \displaystyle a_m \le a_{m-1}\left(1-\frac{D_{\phi}d_m \tilde C}{b_{m-1}}\right) $ \end{center} where we define \begin{center} $ \displaystyle \tilde D_{\phi}:=(1-\phi)(1-\phi/2)\tilde C $ \end{center} so that $\tilde C$ gets absorbed. One may ask if one could proceed also with $D_{\phi}$, but then the third-but-last display on page 424 in \cite{bu} would require \begin{center} $ \displaystyle 1+D_{\phi}d_m/b_{m-1}-D_{\phi}d_m\tilde C/b_{m-1}-D_{\phi}^2d_m^2 \tilde C/b_{m-1}^2 \le 1 $ \end{center} which is not evident (note that the same recursion for the $b_m$ as in \cite[p. 424]{bu} is valid). Maybe one could distinguish between the case where the inequality holds, i.e., where \begin{center} $ \displaystyle \tilde C \ge \frac{1}{1+D_{\phi}d_mb_{m-1}} $ \end{center} holds and the other case, but we do not see any advantage. \ \\

However, we can also apply Temlyakov's lemma \cite[Lem. 12.3]{bu} and get \begin{center} $ \displaystyle a_m^{2+D_{\phi}} \le (1+C_{\phi}^2 m \tilde C)^{-\tilde D_{\phi}} $ \end{center} as analog to the last inequality on page 424 in \cite{bu}. The bound $C_{\phi}^{-2\tilde D_{\phi}} \le 2$ is still valid and we get the inequality \begin{center} $ \displaystyle a_m^{2+\tilde D_{\phi}} \le 2m^{-\tilde D_{\phi}} $ \end{center} on the set $\tilde B_n(m)$. \ \\

On the complement of $\tilde B_n(m)$, we follow exactly the same steps as in \cite{bu}, resulting in the factor $(\tilde C)^{-1}$ before the first summand of the right hand side in the last display on page 425.          \ \\

Finally, we get \vspace{0.15cm} \begin{center} $ \displaystyle ||\hat R^mf||_{(n)}^2 \le \max\left(2m^{-\frac{(1-\phi)(1-\phi/2)\tilde C}{2+(1-\phi)(1-\phi/2)\tilde C}}, 2(\tilde C)^{-1} \phi^{-1} (1-\phi/2)^{-1} \Delta_n (||\beta_n||_1+m\gamma_n)+m\Delta_n \right) $ \end{center} \vspace{0.15cm}
with $\gamma_n$ as in Equation (12.30) in \cite{bu}. \ \\

Using the fact that $\phi$ and $D_{\phi}$ are fixed, that $\gamma_n=\mathcal{O}_P(n^0)$ and $\Delta_n=\mathcal{O}_P\left(\sqrt{\ln(p_n)/n}\right)$ and invoking (P2) and the assumptions on $(m_n)_n$, we asymptotically conclude that \begin{center} $ \displaystyle ||\hat R^mf||_{(n)}^2=o_P(n^0). $ \end{center} \begin{flushright} $_\Box $ \end{flushright} \end{bew}

\begin{rem} In the proof, it might be tempting to absorb $\tilde C$ (or $a$) already in $\phi$ from Lemma \ref{122surro}. However, it turns out that this would not reflect the impact of the constant $\tilde C$ appropriately. \end{rem}

\begin{rem} The last line in the proof of Theorem \ref{asypredsing} means that although the convergence rate is indeed slower for a fixed $n$ due to the constant $\tilde C$ that enters the exponent resp. that enters as factor in the arguments of the maximum, \textbf{we asymptotically do not lose anything compared to $L_2-$Boosting when performing SingBoost}! \end{rem}  

\begin{rem} Again, we did not perform a case-by-case-analysis w.r.t. $m$ using that only each $M-$th iteration is a singular iteration. However, since we get a bound for each $||R^mf||_{(n)}^2$, we do not see any advantage in distinguishing cases w.r.t. $m$ here. \end{rem}   \ \\

\renewcommand\refname{References}
\bibliography{Biblio}
\bibliographystyle{abbrvnat}
\setcitestyle{authoryear,open={((},close={))}}

\end{small}

\end{document}